\documentclass{icm2010}

\usepackage{amsmath}
\usepackage{amssymb}

\usepackage[pdftex,colorlinks,urlcolor=blue,pdfstartview=FitH]{hyperref}

\newcommand{\Rm}{\ensuremath{\mathbb{R}}}

\newcommand{\mC}{\ensuremath{\mathcal{C}}}

\newcommand{\mT}{\ensuremath{\mathcal{T}}}

\newcommand{\mW}{\ensuremath{\mathcal{W}}}

\newcommand{\modo}{\!\!\!\mod}

\newcommand{\Nm}{\ensuremath{\mathbb{N}}}
\newcommand{\Zm}{\ensuremath{\mathbb{Z}}}

\newcommand{\Tm}{\ensuremath{\mathbb{T}}}

\newcommand{\e}{\ensuremath{\epsilon}}

\newtheorem{lem}{Lemma}
\newtheorem{thm}{Theorem}
\newtheorem{conj}[lem]{Conjecture}

\newtheorem{cor}[lem]{Corollary}
\newtheorem{prop}{Proposition}

\newtheorem{defn}[lem]{Definition}

\def\proof {\noindent{\sc{Proof. }}}
\def\qed {\mbox{}\hfill {\small \fbox{}} \\}
\def\lto{\longrightarrow}
\def\lmto{\longmapsto}

\def\leq{\leqslant}
\def\geq{\geqslant}
\title[Arnold diffusion]{Arnold's diffusion: from the \textit{a priori} unstable to the 
\textit{a priori} stable case.}
\author[Patrick Bernard]{Patrick Bernard
\thanks{Membre de l'IUF}%
}
\contact[patrick.bernard@ceremade.dauphine.fr]{CEREMADE, UMR CNRS 7534,
Place du Marechal de Lattre de Tassigny, 75775 Paris cedex 16, France}

\begin{document}

\begin{abstract}
We expose some  selected topics 
concerning the  instability  of the action variables in  \textit{a priori} 
unstable Hamiltonian systems, and outline a new strategy
that may allow to apply these methods to \textit{a priori} stable systems.
\end{abstract}

% Your AMS 200 Classification should come here %
\begin{classification}
37J40, 37J50, 37C29, 37C50, 37J50.
\end{classification}

% Any keywords %
\begin{keywords}
Arnold's diffusion, normally hyperbolic cylinder, partially hyperbolic tori,
homoclinic intersections, Weak KAM solutions, variational methods, action 
minimization.
\end{keywords}

% Do not remove next line %
\maketitle

\section{Introduction}
A very classical problem in dynamics consists in studying
the Hamiltonian system on the symplectic manifolds
 $T^*\Tm^n=\Tm^n\times \Rm^n$
generated by the Hamiltonian
\begin{align}\nonumber
H_{\e}:\Tm\times T^*\Tm^n=\Tm\times\Tm^n \times \Rm^n&\lto \Rm\\
\label{stable}
(t,x,y)&\lmto \frac{1}{2}\|y\|^2+\e G(t,x,y)
\end{align}
where $\e$ is a small perturbation parameter.
More general unperturbed systems $h(y)$ can be considered
instead of  $\|y\|^2/2$, but we restrict to that
particular case in the present paper in order to simplify some notations.
For $\e=0$, the system is integrable, and the momenta $y$
are integrals of motion.
For $\e>0$, these variables undergo small oscillations.
KAM theory implies that these oscillations  remain permanently
bounded  for many initial conditions.
For other initial conditions, a large evolution might be possible. 
By Nekhoroshev theory, it  must  be extremely slow.
The questions we discuss in the present text is whether this large
evolution is actually possible, and to what geometric structures
it is associated.

Let us consider a resonant momentum $y_0=(I_0,0)\in \Rm^m\times \Rm^r=\Rm^n$,
and assume that $I_0$ is not resonant, which means that 
$k\cdot I_0$ never belongs to  $\Zm$ for $k\in \Zm^m, k\neq 0$.
In order to study the dynamics near the torus $\{y=y_0\}$,
it is useful to introduce the notations  
$x=(\theta,q)\in \Tm^m\times \Tm^r$,
and $y=(I,p)\in \Rm^m\times \Rm^r$, $m+r=n$.
In the neighborhood of the torus $\{y=y_0\}$, the dynamics is 
approximated by the averaged system
$$
\frac{1}{2}\|y\|^2+\e  V(q),
$$
where 
$$
V(q)=\int_{\Tm\times \Tm^m} G(t,\theta, q, y_0) d\theta dt.
$$
Following a classical idea of Poincar\'e and Arnold,
we can try to exploit this observation by considering the system
\begin{equation}\label{unstable}
H(t,\theta,q,I,p)=\frac{1}{2}
\|p\|^2+\frac{1}{2}
\|I\|^2 -\e V(q)-\mu R(t,\theta, q,I,p)
\end{equation}
with  a second perturbation parameter $\mu$ independent from $\e$. We assume that $V$ has a unique non-degenerate 
minimum, say at $q=0$.
Fixing $\e >0$, we can  study this system for $\mu>0$ small enough,
which is a simpler problem which may give some hints about the dynamics
of (\ref{stable}).
The reason why instability is more easily proved 
in (\ref{unstable}) than in 
 (\ref{stable})
is the presence of the hyperbolic fixed point at $(0,0)$
of the $(q,p)$ component of the averaged system.
Studying (\ref{unstable}) for $\mu>0$ small enough is thus
called the \textit{a priori} unstable problem, or the 
\textit{a priori} hyperbolic problem. In contrast,
the Hamiltonian (\ref{stable}) is called  \textit{a priori} stable.
The \textit{a priori} unstable case is by now quite well understood  
for $m=1$, see \cite{T, CY1, CY2, DLS, JAMS, X} for example. The 
\textit{a priori} unstable case for $m>1$ and the
\textit{a priori} stable case can be considered as widely
open, in spite of the important announcements of John Mather in \cite{M:A}.
The starting point in the study of (\ref{unstable}) is  the famous 
paper of Arnold, \cite{A}. In  this paper, Arnold introduced a particular
\textit{a priori} unstable system
where some geometric structures associated to diffusion,
partially hyperbolic tori (that he called whiskered), their stable and
unstable manifolds, and heteroclinic connections,  can be almost 
explicitly described.
This geometric structure have been called a transition chain.
Most of the subsequent works on the \textit{a priori} unstable problem  have consisted
in  trying to find transition chains in more general
cases, but understanding the general \textit{a priori} unstable 
Hamiltonian have required a  change of paradigm: 
from partially hyperbolic tori to normally hyperbolic cylinders.   
The variational methods introduced by John  Mather in \cite{Ma:93} and 
Ugo Bessi
in \cite{Bessi} have also  been very influential.

Transforming the understanding gained on the dynamics of 
(\ref{unstable}) to informations on the \textit{a priori} stable 
case is not an easy task.
Since we understand the system  (\ref{unstable}) when $m=1$
the first attempt should be to study (\ref{stable}) in the neighborhood
of an $(n-1)$-resonant line, for example the line  consisting of momenta
of the form $y=(I,0), I\in \Rm$. We could hope to prove the existence
of drift along such a line by using the  \textit{a priori} unstable 
approximations near each value of $y$.
However, we face the problem that an approximation like (\ref{unstable})
 holds only in the neighborhood of the torus 
$\{y=(I,0)\}$ when the frequency
$I\in \Rm$ is irrational.
Near the torus $\{y=(I,0)\}$  with $I$ rational, one should 
use an approximation of the form
$$
H(t,x,y)=\frac{1}{2}\|y\|^2-\e W(x)-\mu R(t,x,y)
$$
and different methods must be used.
This is often called the problem of double resonances when $n=2$.
We will call it the problem of \textit{maximal resonances}.

Our general goal in this paper is to study \textit{a priori}
unstable systems with a sufficient generality to be able to
gain informations on the \textit{a priori} stable case.
We start with a relatively detailed description of the Arnold's
example 
in Section \ref{Arnold}, which is also  an occasion to settle some notations
and introduce some important objects, like the partially hyperbolic tori,
their stable and unstable manifolds, and the associated generating
functions. Working with these generating functions allows to highlight
the  connections between the various classical approaches, geometric 
methods, variational methods, and weak KAM theory.
Then, from the end of Section \ref{Arnold} to Section \ref{general},
we  progressively generalize the setting and indicate how the 
methods introduced on the  example of Arnold can be improved to 
face the new occurring difficulties. We  present the Large Gap Problem,
which prevents Arnold's mechanism from  being directly
 applied to general \textit{a priori} unstable
systems, and explain how the presence of a normally
hyperbolic cylinder can be used to solve this Problem
and prove instability in general \textit{a priori} unstable systems.
In section \ref{hyper} we give   a new result from  \cite{hyperbolic}, on the existence
of normally hyperbolic cylinders in the \textit{a priori}
stable situation, which should allow to apply
the tools exposed in the previous sections
to \textit{a priori} stable systems.
This suggests a possible strategy to prove the following conjecture:

\begin{conj}
For a typical perturbation $G$, there exists two positive  numbers $\epsilon_0$ and $\delta$, such that, for each $\e\in ]0,\e_0[$, 
The system (\ref{stable}) has an orbit 
$$(\theta(t),q(t),\dot \theta(t), \dot q(t)):\Rm\lto 
\Tm\times \Tm^{n-1}\times \Rm
\times \Rm^{n-1}
$$ such that 
$\sup_t \dot \theta -\inf_t \dot \theta >\delta$.
\end{conj}

We are currently working on  this  program in collaboration with Vadim 
Kaloshin and Ke Zhang.
The same conjecture can be stated with a more general unperturbed 
system $h(y)$, and the same proof should work provided $h$ is convex
and smooth.
Our strategy of proof does not consist in solving
the maximal resonance problem, but rather in observing that the 
conjectured statement can be reached without solving
that difficulty. In that respect, what we expose is much easier
than the project of Mather as announced in \cite{M:A}.
The result is  weaker since only limited diffusion is obtained.
The maximal resonance problem has to be solved in order to
prove the existence of global diffusion on a whole 
resonant line, or  from one resonant line
to another.
Our strategy, on the other hand, has the advantage of working 
with all $n\geq 2$, while Mather is limited to $n=2$ at the moment.

\section{The example of Arnold and some extensions}\label{Arnold}
Following Arnold \cite{A}, we consider the Hamiltonian
\begin{equation}\label{ae}
H(t,\theta,q,I,p)=\frac{1}{2}
\|p\|^2+\frac{1}{2}
\|I\|^2 +\e (\cos (2\pi q)-1)(1+\mu f(t,\theta,q))
\end{equation}
with $(t,\theta,q,I,p)\in \Tm\times \Tm\times \Tm\times \Rm\times\Rm$.
We will often use the corresponding Lagrangian 
$$
L(t,\theta,q,\dot \theta, \dot q)=\frac{1}{2}
\|\dot q\|^2+\frac{1}{2}
\|\dot \theta\|^2 +\e (1-\cos (2\pi q))(1+\mu f(t,\theta,q)).
$$
We will see $\e >0$ as a fixed parameter, 
and discuss mainly the small parameter $\mu$.
When $\mu=0$, the variable $I$ is an integral of motion. Our goal 
is to study its evolutions for $\mu>0$.
The form of the perturbation is chosen in such a way that
 the two-dimensional tori
$$
\mT(a)=\Tm^2\times \{0\}\times\{a\}\times \{0\}, \quad a\in \Rm
$$ 
are invariant in the extended phase space,  and carry a linear motion of frequency $(1,a)$.
By studying invariant manifolds attached to these invariant tori,
Arnold discovered  a remarkable  diffusion mechanism, now called the Arnold Mechanism,
that we are  now going to  describe.
In the case $\mu=0$, the tori $\mT(a)$ 
 appear as the products of the hyperbolic fixed
point $\{0,0\}$ of the pendulum in $(q,p)$ by the invariant torus 
$\Tm\times \Tm\times \{a\}$ of the integrable system in the $(t,\theta, I)$ space.
They are thus partially hyperbolic, and have  stable and unstable manifolds,
which coincide and can be given explicitly as
$$
\mW(a)=\big\{ \big(t,\theta, q, a, \pm \partial_q S_0(q)\big): (t, \theta, q)\in \Tm^3 \big\}
$$
with
\begin{equation}\label{S0}
S_0(q)=\frac{2\sqrt{\e}}{\pi}(1-\cos (\pi q)).
\end{equation}
The coincidence and compactness of these
stable manifolds is a very special feature of the unperturbed case $\mu=0$.
For $\mu\neq 0$, the tori $\mT(a)$ still have stable and unstable manifolds 
which can be described as follows:
There exists two functions 
\begin{equation}\label{S}
S^{\pm}_{a,\mu} (t,\theta,q):\Tm\times \Tm\times [-3/4,3/4]\lto \Rm,
\end{equation}
which converge to $\pm S_0$ when $\mu \lto 0$, and
such that the graphs
$$
\mW^{\pm}_{\mu}(a)
=\big\{\big(t,\theta,q \mod 1,a+ \partial_\theta S^{\pm}
(t,\theta, q), \partial_q S^{\pm}(t, \theta, q)\big)\big\}
$$
are pieces of the stable and unstable manifolds of the torus $\mT(a)$.
More precisely, the set 
$\mW^+(a)$ is negatively invariant under the extended Hamiltonian flow,
and 
$$\mT(a)=\bigcap_{t\leq 0} \varphi^t\big( \mW^+(a))
$$
while 
the set 
$\mW^-(a)$ is positively invariant under the extended Hamiltonian flow,
and 
$$\mT(a)=\bigcap_{t\geq 0} \varphi^t\big( \mW^-(a)).
$$
The functions $S^{\pm}_a$ solve the Hamilton-Jacobi equation
$$
\partial_t S+H(t,\theta, q,a+\partial_\theta S,\partial_q S)=a^2/2,
$$ 
which merely says that the invariant manifolds are contained in the energy
level of the torus.
The functions $S^{\pm}_a$ have an expression in terms of the action:
\begin{align}\label{Saction}
S^+_a(t,\theta,q)&=
\int_{-\infty}^{\tau}L(s,\theta^+(s),q^+(s),\dot \theta^+(s),\dot q^+(s))
-a\dot \theta^+(s) +a^2/2ds\\
\nonumber
S^-_a(t,\theta,q)&=
\int_{\tau}^{+\infty}L(s,\theta^-(s),q^-(s),\dot \theta^-(s),\dot q^-(s))
-a\dot \theta^-(s)+a^2/2ds,
\end{align}
where $\tau$ is any real number such that $\tau \mod 1=t$, and 
$(\theta^{\pm}(s),q^{\pm}(s))$ is the solution of the Euler-Lagrange equations
such that 
$$\theta^{\pm}(\tau)=\theta, q^{\pm}(\tau)=q\mod 1, \dot \theta^{\pm}(\tau)=a+\partial_{\theta} S^{\pm}(t,\theta,q),
\dot q^{\pm}(\tau)=\partial_{q}S^{\pm}(t,\theta,q).
$$
Note that the result does not depend on the choice of $\tau$.

\subsection{Homoclinic orbits}\label{hom}
If $(T,\Theta,Q)\in \Tm\times \Tm \times ]1/4,3/4[$
is a critical point of the function
$$
\Delta_a (t,\theta,q)=S^+_a(t,\theta,q)-S^-_a(t,\theta, q-1),
$$
then the point 
\begin{align*}
&\big(T,\Theta, Q\mod 1,a+\partial_{\theta}S_a^+(T,\Theta, Q),
\partial_qS_a^+(T,\Theta, Q)\big)\\
=&
(T,\Theta, (Q-1)\mod 1,a+\partial_{\theta}S_a^-(T,\Theta, Q-1),
\partial_qS_a^-(T,\Theta, Q-1)\big)
\end{align*}
obviously belongs both to $\mW^+(a)$ and $\mW^-(a)$, hence it is a homoclinic point.
It is a transversal homoclinic point if in addition the Hessian
of $\Delta_a$ has rank two (it can not have rank 3 because the intersection is necessarily one-dimensional).
It is not obvious at this point that the function $\Delta_a$
necessarily has critical points on the domain $\Tm\times \Tm \times ]1/4,3/4[$.
When $\mu$ is small enough,  this follows from:

\begin{lem}
If $(T,Q)$ is a critical point of the function 
 $\bar \Delta_a:(t,q)\lmto\Delta_a(t,q,1/2)$, then  $(T,Q, 1/2)$ is a critical point 
of $\Delta_a$, 
 hence the manifolds $\mW^-(a)$ and $\mW^+(a)$ intersect
above $(T,\Theta ,1/2)\in \Tm^3$.  This homoclinic point 
is transversal if and only if the Hessian of $\bar \Delta_a$ at $(T,Q)$
is a non-degenerate $2\times 2$ matrix.
\end{lem}

Note that the function $\bar \Delta_a $ is defined on $\Tm^2$,
and therefore it has critical points.

\proof
We have $\partial_tS^+(T,Q,1/2)=\partial_tS^-(T,Q,-1/2)$,
let us denote by $e$ this value.
We also have  $\partial_{\theta}S^+(T,Q,1/2)=\partial_{\theta}S^-(T,Q,-1/2)$,
we denote by $I$ this value.
It is enough to prove that $\partial_qS^+(T,Q,1/2)=\partial_qS^-(T,Q,-1/2)$.
In order to do so, it is enough to observe that $\partial_qS^+$
is the only non-negative solution of the equation
$$e+H(T,\Theta, 1/2,a+ I,.)=a^2/2,$$
and that precisely the same characterization is true for 
$\partial _qS^-(T,Q,-1/2)$.
Note that the equation above has two solutions, and that we can discriminate 
between them because we work in a perturbative setting which gives 
us rough informations on the signs. In more general situation, this
is a source of difficulty.
\qed

\subsection{Heteroclinic orbits}\label{hetero}
We have proved the existence of homoclinic orbits. But what
is interesting for Arnold diffusion are heteroclinic orbits between
different tori.
We can deduce the existence of a heteroclinic orbit
between $\mT(a)$ and $\mT(a')$ provided we can find a critical point 
of the function 
$$\Tm\times \Rm\times ]1/4,3/4[\ni (t,\theta, q)\lmto
S_a^{+}(t,\theta,q)-S_{a'}^-(t,\theta, q-1)+(a-a')\theta,
$$
where we have lifted the functions $S$ without changing their names.
As before, we can limit ourselves to finding critical points of the function
\begin{equation}\label{Sigma}
\Sigma_{a,a'}:\Tm\times \Rm\ni (t,\theta)\lmto
S_a^{+}(t,\theta,1/2)-S_{a'}^-(t,\theta, -1/2)+(a-a')\theta,
\end{equation}
but the term $(a-a')\theta$ prevents us from finding them 
using  a  global variational method when $a'\neq a$.
This  reflects the fact that we are studying a non exact
Lagrangian intersection problem.
For $\mu=0$, heteroclinics do not exist.
However, recalling that 
$\bar \Delta_a(t,\theta)=\Delta_a(t,\theta,1/2)$, we have:

\begin{lem}\label{heteroclinics}
If the function $\bar \Delta_a(t,q)$ has a non-degenerate critical point,
then the functions $\Sigma_{a,a'}$ and $\Sigma_{a',a}$ both
have a non-degenerate critical point provided $a'$ is sufficiently close to $a$.
\end{lem}

\proof
The theory of partial hyperbolicity implies that the stable and unstable
manifolds $\mW^{\pm}_{\mu}(a)$ depend regularly on the parameter $a$.
As a consequence, their generating functions $S^{\pm}_a$
also regularly depend on $a$, and the functions 
 $\Sigma_{a,a'}$ depend regularly on $a$ and $a'$.
The result follows since $\Sigma_{a,a}=\bar \Delta_a$.
\qed

We say that $a_0, a_1,\ldots, a_k$ is an elementary transition chain
if the functions $\Sigma_{a_{i-1},a_i}$ have non-degenerate critical points.
We will sometimes use the same terminology for the different requirement that 
these functions have isolated local minima.
From Lemma \ref{heteroclinics}, we deduce:

\begin{prop}\label{chain}
Let $\mu$ be given and sufficiently small.
Let $[a^-, a^+]$ be an interval such that each of the functions 
$\bar \Delta_{a,\mu}, a\in [a^-, a^+]$ have a non-degenerate critical point, which means
that each of the tori $\mT_{\mu}(a),a\in[a^-, a^+]$ has a transversal homoclinic
orbit. Then there exists an elementary transition chain 
$a^-=a_0,a_1, \ldots , a_k=a^+$.
\end{prop}

\proof
Let us consider the set $A\subset [a^-,a^+]$ of points that can be reached from $a^-$
by a transition chain.
The set $A$ is open : If $a'\in A$, then there exists a transition chain
$a^-=a_0, a_1,\ldots a_k=a'$ and, by Lemma \ref{heteroclinics},
the sequence  $a^-=a_0, a_1,\ldots a_k, a_{k+1}$
is a transition chain when $a_{k+1}$ is sufficiently close to $a$.
The set $A$ is closed : Let $a$ be  in the closure  of $A$. By Lemma 
\ref{heteroclinics}, the pair $a, a'$ is a transition chain when $a'$
is close to $a$. Since $a$ is in the closure of $A$, the point $a'$ can be taken in $A$. Then,  there exists a transition chain
$a^-=a_0,\ldots, a_k=a'$, and then the longer sequence
$a^-=a_0,\ldots, a_k, a_{k+1}=a$ is a transition chain between $a_0$
and $a$, hence $a\in A$.
Being open, closed and not empty (it contains $a_0$), the set $A$ is equal
to $[a^-,a^+]$.
\qed

The existence of transition chains implies the existence 
of diffusion orbits.
This is proved by Arnold invoking an ``obstruction property''.
This obstruction property is a characteristic of the local dynamics
near the partially hyperbolic tori.
It has been proved by Jean-Pierre Marco in \cite{Ma:96},
see also \cite{Cr:97,FM:00}.
The most appealing way to understand the geometric shadowing 
of transition chains is to use the following statement of 
Jacky Cresson \cite{Cr:00}, which  can be seen as a strong obstruction
property:

\begin{lem}\label{transitivity}
If there exists a transversal heteroclinic between $\mT(a)$ and $\mT(a')$
and a transversal heteroclinic between $\mT(a')$ and $\mT(a'')$, then there
exists a transversal heteroclinic between $\mT(a)$ and $\mT(a'')$.
\end{lem}
This Lemma implies:

\begin{cor}
If $a_0, a_1, \ldots, a_k$ is an elementary transition chain, then there exists 
a transversal heteroclinic orbit between $\mT(a_0)$ and $\mT(a_k)$.
\end{cor}

Putting everything  together, we obtain:
\begin{thm}\label{t1}
Let $\mu$ be given and sufficiently small.
Let $[a^-, a^+]$ be an interval such that each of the functions 
$\bar \Delta_{a,\mu}, a\in [a^-, a^+]$ have a non-degenerate critical point.
Then there exists a heteroclinic orbit between $\mT(a^-)$ and $\mT(a^+)$.
\end{thm}

\subsection{Poincar\'e-Melnikov approximation}\label{PM}
We have constructed diffusion orbits under the assumption that 
transversal homoclinics  exist.
We have proved that homoclinic orbits necessarily exist,
and one may argue that transversality should hold for typical systems,
we will come back on this later.
However, it is useful to be able to check whether transversality holds
in a given system.
A classical approach consists in proving the existence
of non-degenerate critical points of the functions
$\bar \Delta_{a,\mu}$ defined in Lemma \ref{heteroclinics} by  expanding them
in power series of $\mu$.
As a starting point the generating functions $S^{\pm}_a$ can be expanded as follows:
\begin{align}\label{power}
S^+_a(t,\theta,q)&=S_0(q)+\mu M^+_a(t,\theta,q)+O(\mu^2),\\
\nonumber
S^-_a(t,\theta,q)&=-S_0(q)-\mu M^-_a(t,\theta,q)+O(\mu^2),
\end{align}
where $S_0(q)=\frac{2\sqrt{\e}}{\pi}(1-\cos (\pi q))$ is the generating function
of the unperturbed manifolds, and $M^{\pm}$ are the so-called 
Poincar\'e-Melnikov integrals, 
\begin{align*}
M^+_a(t,\theta,q)&=\e\int_{-\infty}^{t}
F\Big(s,\theta+a(s-t), \frac{2}{\pi}\arctan\big(
e^{2\pi\sqrt{\e}(s-t)}\tan (\pi q/2)\big)\Big) ds
\\
M^-_a(t,\theta,q)&=
\e \int_{t}^{+\infty}
F\Big(s,\theta+a(s-t), \frac{2}{\pi}\arctan\big(
e^{2\pi\sqrt{\e}(t-s)}\tan (\pi q/2)\big)\Big) ds
\end{align*}
where $F(t,\theta, q)=(1-\cos (2\pi q))f(t,\theta, q)$.
To better understand these formula, it is worth
recalling that 
$$s\lmto \frac{2}{\pi} \arctan\big(
e^{2\pi\sqrt{\e}(s-t)}\tan (\pi q/2)\big) $$
is the homoclinic orbit of the 
system $\|p\|^2/2+\e (\cos (2\pi q)-1)$
which takes the value $q$ at time $t$.
The  formula above  are similar to (\ref{Saction}), but the integration is
performed on unperturbed trajectories, which are explicitly known.
For $q\in ]1/4,3/4[$, we obtain:
$$
\Delta_a(t,\theta,q)=S^+_a(t,\theta,q)-S^-_a(t, \theta, q-1)
=\mu M_a(t,\theta, q)+O(\mu^2),
$$
where $M_a$ is the Poincar\'e-Melnikov integral
\begin{align*}
M_a(t,\theta, q)&=M^+_a(t,\theta,q)+M^-_a(t,\theta, q-1)\\
&=\e \int_{\Rm} F\Big(s,\theta+a(s-t), \frac{2}{\pi}\arctan\big(
e^{2\pi\sqrt{\e}(s-t)}\tan (\pi q/2)\big)\Big) ds.
\end{align*}
In the specific case studied by Arnold, where 
$f(t,\theta, q)=\cos (2\pi \theta)  +\cos (2\pi t)$, the Melnikov integral can be computed explicitly through residues, we obtain:
$$
M_a(t,q,1/2)=
\frac{a}{\text{sh}(\pi a/2\sqrt{\e})}\cos (2\pi \theta)+
\frac{1}{\text{sh}(\pi /2\sqrt{\e})}\cos (2\pi t),
$$
it has a non-degenerate minimum at $(t,q)=(0,0)$.
We can conclude, following Arnold:
\begin{thm}[Arnold, \cite{A}]\label{at}
Let us consider the Hamiltonian (\ref{ae}) with
 $f(t,\theta, q)=\cos (2\pi \theta)+\cos (2\pi t)$ and $\mu>0$ small enough.
Given two real numbers $a^-<a^+$, there exists an orbit
$(\theta(t), q(t),I(t),p(t))$ and a time $T>0$ such that 
$I(0)\leq a^-$ and $I(T)\geq a^+$.
\end{thm}

\subsection{Bessi's variational mechanism}\label{var}
Ugo Bessi introduced in \cite{Bessi} a very interesting  approach to study 
the system (\ref{ae}), see also \cite{Bessi:97, Bessi:double}.
 In order to describe this approach, let
 us define the function 
\begin{align*}
A_a&:\Rm\times \Tm\times ]1/4,3/4[\times \Rm\times \Tm\times ]1/4,3/4[\lto \Rm\\
&((t_1,\theta_1,q_1),(t_2,\theta_2,q_2))\lmto
\min \int_{t_1}^{t_2} L(s, \theta(s), q(s),\dot \theta (s), \dot q(s))-a\dot \theta(s)+a^2/2\,ds,
\end{align*}
where the minimum is taken on the set of $C^1$ curves 
$(\theta(s),q(s)):[t_1,t_2]\lto \Tm\times \Rm$ such that 
$$(\theta(t_1),q(t_1))=(\theta_1, q_1-1) \quad \text{and} \quad
(\theta(t_2),q(t_2))=(\theta_2, q_2).
$$
When the time interval $t_2-t_1$ is very large, the minimizing
trajectory in the definition of $A_a$  roughly looks like the concatenation of 
an orbit  positively asymptotic to $\mT(a)$ followed by
  an orbit negatively asymptotic to $\mT(a)$.
Using this observation, and recalling the formula (\ref {Saction}),
it is  possible to prove that
\begin{align*}
&A_a((t_1,\theta_1,q_1),(t_2+k,\theta_2,q_2))\lto\\
&S^+_a(t_2\modo 1, \theta_2 , q_2)-S^-_a(t_1\modo 1, \theta_1, q_1-1)
\end{align*}
when $k\lto \infty$.
Fixing the real numbers $a_0,a_1, \ldots, a_k$ and the integers 
$\tau_1,\ldots ,\tau_k$, we consider
 the discrete action functional 
\begin{align*}
&S^+_{a_0}(t_1\modo 1,\theta_1 \modo 1 ,1/2)+
(a_0-a_1)\theta_1\\
+ &A_{a_1}((t_1,\theta_1  \modo 1, 1/2),(t_2+\tau_2,\theta_2 \modo 1,1/2))
+(a_1-a_2)\theta_2\\
+&A_{a_2}((t_2,\theta_2 \modo 1, 1/2),(t_3+\tau_3,\theta_3 \modo 1,1/2))+(a_2-a_3)\theta_3\\
+&\cdots\\
+&A_{a_{k-1}}(t_{k-1}+\tau_{k-1},\theta_{k-1} \modo 1, 1/2),
(t_{k},\theta_{k} \modo 1,1/2))+(a_{k-1}-a_k)\theta_k\\
-& S^-_{a_k}(t_k\modo 1,\theta_{k} \modo 1,1/2)
\end{align*}
defined on $(]-1,1[\times ]-1,1[)^{k}$.
It is not hard to check that local minima of this discrete action functional
give heteroclinics between the Torus $\mT(a_0)$ and the torus $\mT(a_k)$.
In order to prove that local minima exist, observe that 
 this functional
is approximated by 
$$
\Sigma_{a_0,a_1}(t_1 \modo 1,\theta_1 )+\cdots +
\Sigma_{a_{k-1},a_k}(t_k \modo 1, \theta_k)
$$
%\begin{align*}
%%(a_1-a_0)\theta_1+\\
%&S^+_{a_1}(t_2 \modo 1,\theta_2 \modo 1,1/2)-S^-_{a_2}(t_3\modo 1,\theta_3 \modo 1,1/2)
%+(a_2-a_1)\theta_2+\\
%&\cdots\\
%&S^+_{a_{k-1}}((t_{k}\modo 1,\theta_{k}\modo 1,1/2)-
%S^-_{a_{k}}(t_{k} \modo 1,\theta_{k} \modo 1, 1/2)+(a_k-a_ {k-1})\theta_k.
%\end{align*}
when the integers $\tau_i$ are large enough, with the functions $\Sigma$ as defined in (\ref{Sigma}).
This limit functional has the remarkable structure that the variables 
$(t_i,\theta_i)$ are separated.
This break-down of the action functional into a sum of independent functions
is sometimes called an anti-integrable limit, it is related to 
the obstruction property of the invariant tori, to the $\lambda$-Lemma,
and to the Shilnokov's Lemma, see \cite{Bo:00}.
The limit functional  has an isolated local minimum  provided each of the functions 
$\Sigma_{a_{i-1},a_i}$ has one, which is equivalent to say that 
$a_0, a_1, \ldots, a_k$ is an elementary transition chain.
In this case, the integers $\tau_i$ can be chosen large enough so that 
the action functional above has a local minimum, which gives a heteroclinic 
orbit between $\mT(a_0)$ and $\mT(a_k)$.
Technically, this method has several advantages. 
In our presentation
we introduced the generating functions $S^{\pm}_ {a_i}$ of the invariant
manifolds of the involved tori in order to stress the relations between 
the geometric and the variational method, and also
because this works in a more general setting, see \cite{Bo:00}
for example.
In our context and when $\mu>0$ is small enough,
 it is easier to directly approximate 
the functions $A_a$ in terms of the Melnikov integrals,
and to use the following approximation for the action functional with large $\tau_i$
and small $\mu$
without the intermediate step through $S^{\pm}$:
\begin{align*}
\mu M_{a_0}(t_0 \modo 1, \theta_0\modo 1, 1/2)+(a_1-a_0)\theta _0+\cdots+\\
\mu M_{a_k}(t_k\modo 1,\theta_k\modo 1,1/2)+(a_k-a_{k-1})\theta_k.
\end{align*}
The corresponding calculations, performed in \cite{Bessi}, are much more elementary
than those required to derive the expansions (\ref{power}).

\subsection{Remarks on estimates} We have up to that 
point carefully avoided to discuss the subtle and important aspect 
of explicit estimates. 
In order to complete rigorously the proof
of Theorem \ref{at}, we should prove the existence of a threshold
$\mu_0(\e)$ such that the Melnikov approximation holds,
simultaneously for all $a$, when $0<\mu<\mu_0(\e)$.
This can actually been done, with 
$$\mu_0(\e)= e^{-\frac{C}{\sqrt{\e}}},
$$
but it is not simple,
since it requires to study carefully the expansions of the functions
$S^{\pm}_a$ and how the coefficients depend on $a$ and $\epsilon$.
This is related to the so-called splitting problem,
see \cite{LMS}.
 As we mentioned above
the approach of Bessi allows to prove  that Theorem \ref{at}
holds for $0<\mu<\mu_0(\e)$ without estimating the splitting.

It is also  important to give time estimates, that is to estimate
the time needed for the variable $I$ to perform
a large evolution. Once again, this is closely related to the 
splitting estimates, although these can be avoided by using the 
method of Bessi.
One should distinguish two different problems. Either we fix $\e$,
and try to estimate the time as a function of $\mu$, or we take 
$\mu$ as a function of $\e$, say $\mu=\mu_0(\e)/2$, and
try to estimate the time as a function of $\e$.

The second problem is especially important, because it is relevant for
the study of the \textit{a priori} stable problem. 
Once again, Ugo Bessi obtained the first estimate, 
$$
T=e^{\frac{C}{\sqrt{\e}}}.
$$
Estimating the time on examples allows to test the optimality of 
Nekhoroshev exponents, see \cite{MS,LM, Z} for works in that direction,
see also \cite{BK} concerning the question of time estimates.

It is worth mentioning also that in the first problem, estimating 
the time as a function of $\mu$, the estimate is polynomial, and
not exponentially small. This was first understood by Pierre Lochak, and
 proved by Bessi's method  in \cite{CRAS}, where the estimate 
$T=C/\mu^2$ is given, see also \cite{Cr:01}. The optimal estimate is $T=C|\ln \mu|/ \mu$,
as was conjectured by Lochak in \cite{L} and proved by Berti, Biasco and Bolle
in \cite{BBB}, see also \cite{BB}.

Returning to the question of the threshold of validity,
let us discuss what happens when $\mu$ is increased above $\mu_0(\e)$.
The content of Section \ref{PM} on finding transversal homoclinics
via the Poincar\'e-Melnikov approximation breaks down, but the 
geometric constructions of the earlier sections is still valid.
 Theorem \ref{t1} holds as long as the invariant tori $\mT(a)$
remain partially hyperbolic, and that their stable and unstable manifold
can be represented by generating functions like (\ref{S}).
Actually, the methods we are now going to expose allow even to relax
this last assumption. 
Being able to treat larger values of $\mu$ is especially important in view
of the possible applicability to the \textit{a priori} stable problem.

\subsection{Higher dimensions}\label{hd}
Let us now discuss the following immediate generalization in higher
dimensions of Arnold's example:
$$
H(t,\theta,q,I,p)=\frac{1}{2}
\|p\|^2+\frac{1}{2}
\|I\|^2 -\e V(q)(1+\mu f(t,\theta,q))
$$
with $(t,\theta,q,I,p)\in \Tm\times \Tm^m\times \Tm^r\times \Rm^m\times\Rm^r$,
where $V(q)$ is a non-negative function having a unique non-degenerate minimum
at $q=0$, with $V(0)=0$.
The main difference with the example of Arnold 
appears for $r>1$. In this case,  the system
is not integrable even for $\mu=0$.
There still exists a family  of partially hyperbolic 
tori of dimension $m$,
$$\mT(a):=\{(t,\theta,0,a,0), (t,\theta)\in \Tm\times \Tm^m \}
$$
parametrized by $a\in \Rm^m$, but
 the 
system $\|p\|^2/2- \e V(q)$ is not necessarily integrable any more.
As a consequence we do not know explicitly the stable and unstable manifolds of the 
hyperbolic fixed point $(0,0)$, and so we do not have a  perturbative 
setting to describe the stable and unstable manifolds 
of the hyperbolic tori $\mT(a)$. This is also what happens for $r=1$
if $\mu$ is not small enough.
There is no obvious generalization of the generating
functions $S^{\pm}_a$
in that setting, because the stable and unstable manifolds are not necessarily graphs
over a  prescribed domain.
The proof of the existence of homoclinic orbits 
as given in \ref{hom}  thus breaks down.
The existence of homoclinic orbits in that setting can still be
proved by global variational methods, as is now quite well
understood, see \cite{Bolotin, F:hetero,CP, etds, E:94} for example.

The proof is quite easy in our context, let us give a rapid sketch.
We first define a function  $A_a$  similar to the one appearing in Section \ref{var},
but slightly different:
\begin{align*}
A_a&:\Rm\times \Tm^m\times \Tm^r\times \Rm\times \Tm^m\times \Tm^r\lto \Rm\\
&((t_1,\theta_1,q_1),(t_2,\theta_2,q_2))\lmto
\min \int_{t_1}^{t_2} L(s, \theta(s), q(s),\dot \theta (s), \dot q(s))-a\dot \theta(s)+a^2/2\,ds,
\end{align*}
where the minimum is taken on the set of curves 
$(\theta(s),q(s)):[t_1,t_2]\lto \Tm^m\times \Tm^r$ such that 
$(\theta(t_i),q(t_i))=(\theta_i,q_i)$ for $i=1$ or $2$.
Let us set 
$$
\xi(a):= \liminf _{\Nm\ni k\lto \infty}A_a((0,0,q_0),(k,0,q_1)),
$$
and consider  a sequence of minimizing extremals
$$(\theta_i(t),q_i(t)):[0,k_i]\lto \Tm^m\times \Tm^r$$
 such that 
$(\theta_i(0),q_i(0))=(0,q_0)$,  $(\theta_i(k_i),q_i(k_i))=(0,q_1)$,
$k_i\lto \infty$, 
and 
$$
\int_{0}^{k_i} L(s, \theta_i(s), q_i(s),\dot \theta_i (s), \dot q_i(s))-a\dot \theta_i(s)+a^2/2\,ds\lto \xi(a).
$$
Let $M$ be a submanifold of $\Tm\times \Tm^m\times \Tm^r$ which separates
$\Tm\times \Tm^m\times \{q_0\}$ from $\Tm\times \Tm^m \times \{q_1\}$,
and let $T_i\in[0,k_i]$ be a time such that 
$(T_i\mod 1, \theta_i(T_i), q_i(T_i))\in M$,
and let $\tau_i$ be the integer  part of $T_i$.
It is not hard to check that the curves 
$(\theta_i(t-\tau_i),q_i(t-\tau_i))$
converge (up to a subsequence) uniformly on compact sets to a limit
$(\theta_{\infty}(t),q_{\infty}(t)):\Rm\lto \Tm^m\times \Tm^r$.
This limit curve satisfies 
\begin{equation}\label{xi}
\int_{-\infty}^{\infty}  L(s, \theta_{\infty}(s), q_{\infty}(s),
\dot \theta_{\infty}(s), \dot q_{\infty}(s))-a\dot \theta_{\infty}(s)+a^2/2\,ds
=\xi(a),
\end{equation}
and the corresponding orbit is a heteroclinic from $\mT_0(a)$ to $\mT_1(a)$.
We call minimizing heteroclinics (for the lifted system)
those which have minimal action, or in 
other words those which satisfy (\ref{xi}).
In the original system (before taking the covering),
we call minimizing homoclinic orbit a homoclinic which lifts 
to a minimizing heteroclinic.

Let us now try to establish some connections 
between the present discussion and the proof
of the existence of homoclinic orbits given in \ref{hom}.
We  define two functions on $\Tm\times \Tm^m\times \Tm^r$:
\begin{align}\label{Sa}
S^-_a(t,\theta,q)&=-\liminf_{\Nm\ni k\lto \infty}
A_a\big((t,\theta,q),(k,0,q_1)\big)\\
S^+_a(t,\theta, q)&=
\liminf_{\Nm\ni k\lto \infty}
A_a\big((0,0,q_0),(t+k,\theta,q)\big).
\end{align}
Note that 
$$
S_a^+(0,0,q_1)=-S_a^-(0,0,q_0)=\xi(a).
$$
The functions $S_a^{\pm}$, whose definition
is basic both in Mather's (\cite{Ma:93}) and in Fathi's (\cite{Fathi}) theory,
 share many features with those introduced in (\ref{S}), that's why we use the same name.
Let us state some of  their properties:

The function $S^-_a$ is non-positive and  it vanishes on $\Tm\times\Tm^m\times \{q_1\}$
(and only there). Moreover, it is smooth around this manifold,
which is a transversally non-degenerate  critical manifold. 
Let us chose a small $\delta>0$. The set 
$$
\mW^-_{loc}(a):=\Big\{
\big(t,\theta, q,a+ \partial_{\theta} S^-_a,\partial_q S^-_a(t,\theta,q)
\big),\quad 
S^-_a(t,\theta, q)> -\delta
\Big\}$$
is a positively invariant local stable manifold of $\mT_{1}(a)$.

Similarly, $S^+_a$ is non-negative, it is null on $\Tm\times\Tm^m\times \{q_0\}$,
and smooth around it, and this critical manifold is transversally non-degenerate.
The set 
$$
\mW^+_{loc}(a):=\Big\{
\big(t,\theta, q,a+ \partial_{\theta} S^+_a,\partial_q S^+_a(t,\theta,q)
\big),\quad 
S^+_a(t,\theta, q)<\delta
\Big\}$$
is a negatively invariant local unstable manifold of $\mT_{0}(a)$.

The functions $S^{\pm}_a$ also have a global meaning.
Let us give the details for $S^+$.
For each point $(T,\Theta, Q)$, there exists a real number $\tau\in \Rm$
and at least one solution 
$(\theta(s), q(s)):(-\infty, \tau]\lto \Tm^ m\times \Tm^r$ of the Euler-Lagrange equations 
such that $(\tau \modo 1, \theta(\tau),q(\tau))=(T,\Theta,Q)$, and which is calibrated
by $S^+_a$ in the following sense: The relation 
\begin{align*}
&S^+_a(t \modo 1,\theta(t), q(t))-
S^+_a(s \modo 1,\theta(s), q(s))\\
=&\int_s^t L(\sigma, \theta(\sigma), q(\sigma), \dot \theta(\sigma),\dot  q(\sigma))
-a\dot \theta(\sigma) +a^2/2\,  d\sigma
\end{align*}
holds for all $s<t\leq \tau$. The corresponding orbit is asymptotic  either to $\mT_0(a)$
or to $\mT_1(a)$   when $s\lto -\infty$.
It is not easy in general to determine whether the asymptotic torus 
is $\mT_0(a)$ or $\mT_1(a)$  but the following Lemma is not hard to prove:

\begin{lem}
If $S_a^+(T,\Theta, Q)<\xi(a)$, 
then each calibrated curve 
$$(\theta(s), q(s)):(-\infty, \tau]\lto \Tm^ m\times \Tm^r
$$
satisfying  $(\tau \modo 1, \theta(\tau),q(\tau))=(T,\Theta,Q)$,
is $\alpha$-asymptotic to $\mT_0(a)$,
and satisfies 
$$
\int_{-\infty}^{\tau}
L(\sigma, \theta(\sigma), q(\sigma), \dot \theta(\sigma),\dot  q(\sigma))
-a\dot \theta(\sigma) +a^2/2\,  d\sigma
=S_a^+(T,\Theta,Q).
$$
\end{lem}

If the function $S^+_a$ is differentiable at $(T,\Theta,Q)$ then there 
is one and only one calibrated curve as above, it is characterized by the
equations
$$
\dot\theta(\tau)=a+\partial_{\theta}S^+_a(\tau, \Theta, Q)
,\quad
\dot q(\tau)=\partial_{q}S^+_a(\tau, \Theta, Q).
$$
Formally, the critical points of the difference $S_a^+-S_a^-$ correspond to 
heteroclinic orbits (in the lifted system). 
By studying a bit more carefully the relations between the calibrated curves
and the differentiability properties of the functions $S_a^{\pm}$
(which is one of the central aspects of Fathi's Weak KAM theory, see \cite{Fathi}), 
this idea can be made rigorous as follows:

\begin{lem}
If $(T,\Theta, Q)$ is a local minimum  of the function  $S^+_a-S^-_a$, then both $S^+_a$ and 
$S^-_a$ are differentiable at the point $(T,\Theta,Q)$,
we have
$$
\big(T,\Theta,Q,a+\partial_{\theta} S^-,\partial_q S^-\big)=
\big(T,\Theta,Q,a+\partial_{\theta} S^+,\partial_q S^+\big),
$$
and the orbit of this point is   either a heteroclinic between $\mT_0(a)$ and $\mT_1(a)$
or a homoclinic to $\mT_0(a)$ or to $\mT_1(a)$
in the system lifted to the covering, and thus it projects
to an orbit homoclinic to $\mT(a)$ in the original system.
\end{lem}

Although it  is not obvious \textit{a priori} that a 
local minimum of the function $S^+_a-S^-_a$
exists away from $q=q_0$ and $q=q_1$, this follows from the
existence of minimizing heteroclinics, that we already proved.
 More precisely, we have:
\begin{itemize}
\item The minimal value of $S_a^+-S_a^-$ is $\xi(a)$.
\item The point $(T,\Theta,Q)$ is a global minimum of $S_a^+-S_a^-$ if and only
if either $Q\in \{q_0,q_1\}$ or the orbit of the point 
$
\big(T,\Theta,Q,a+\partial_{\theta} S^-,\partial_q S^-\big)=
\big(T,\Theta,Q,a+\partial_{\theta} S^+,\partial_q S^+\big),
$
is a minimizing heteroclinic between $\mT_0(a)$ and $\mT_1(a)$.
\item The set of minima of the function $S_a^+-S_a^-$ properly
contains $\Tm\times \Tm^m\times \{q_0\}\cup \Tm\times \Tm^m\times \{q_1\}$.
\end{itemize}

As a consequence, the trajectory $(\theta(t),q(t),\dot \theta(t), \dot q(t))$
is a minimizing heteroclinic if and only if 
$(S_a^+-S_a^-)(t \modo 1,\theta(t),q(t))=\xi(a)$ for each $t\in \Rm$
(and if $q(t)$ is not identically $q_0$ or $q_1$).
This minimizing heteroclinic is called isolated
if, for some $t\in \Rm$, the point $(\theta(t),q(t))$
is an isolated minimum of the function 
$$
(\theta,q)\lmto (S_a^+-S_a^-)(t\modo 1, \theta, q).
$$

Now we have proved that the stable and unstable manifolds of the torus
$\mT(a)$ necessarily intersect, let us suppose that there exists a 
compact and connected set $A\subset \Rm^m$
such that the intersection is transversal for $a\in A$.
By a continuity argument as in Proposition \ref{chain},
we conclude  that any two points $a^-$ and $a^+$ in $A$ can be connected
by a transition chain, that is a sequence 
$a_0=a^-,a_1,\ldots,a_n=a^+$ such that the unstable manifold of 
$\mT(a_{i-1})$ transversally intersects the stable manifold of 
$\mT(a_i)$.
We would like to deduce the existence of a transversal heteroclinic orbit
 between $\mT(a^-)$ and $\mT(a^+)$, but I do not know whether the higher
codimensional analog of Cresson's transitivity Lemma \ref{transitivity}
 holds. 
However, the weaker obstruction property proved in \cite{Cr:97, FM:00}
is enough to  imply  the existence of  orbits
 connecting  any neighborhood of $\mT(a^-)$ to any neighborhood of 
$\mT(a^+)$. It is also possible  to 
build shadowing orbits using a variational approach.
We need  the slightly different  assumption
 that $A\subset \Rm^m$ is a compact connected
set such that, for all $a\in A$, 
all the minimizing homoclinics of $\mT(a)$
are isolated.
For each $a^-$ and  $a^+$ in $A$, it is then possible to construct 
by a variational method similar to Section \ref{var}
a heteroclinic orbit between $\mT(a^-)$ and $\mT(a^+)$.

\section{The general  \textit{a priori} unstable case}\label{general}
A very specific feature of all the examples studied so far is that the perturbation
preserves the partially hyperbolic invariant tori $\mT(a), a\in \Rm^m$.
We now discuss the general \textit{a priori} unstable system 
(\ref{unstable}).
\subsection{The Large Gap Problem}
Let us assume that $r=1$ and try to apply the method of 
Section \ref{Arnold}. 
There is no explicit invariant torus any more,
but  KAM methods can be applied to prove the existence of many
partially hyperbolic tori.
More precisely, there exists a diffeomorphism 
$$\omega_{\mu}(a):\Rm^m\lto \Rm^m,
$$
close to the identity, such that an invariant quasiperiodic Torus $\mT_{\mu}(a)$ 
of frequency $\omega_{\mu}(a)$ exists, and is close to $\mT(a)$, provided the frequency
$\omega_{\mu}(a)$ satisfies some Diophantine condition.
Moreover,  for such values of $a$, the local stable and unstable
manifolds $\mW^{\pm}_{\mu}(a)$ can be generated by functions
$$S^{\pm}_{a,\mu} (t,\theta,q):\Tm\times \Tm^m\times [-3/4,3/4]\lto \Rm,
$$
as earlier.
So we have exactly the same picture as in Section \ref{Arnold}, except
that the objects are defined only on a subset $A_{\mu}\subset \Rm^m$
of parameters.
In order  to reproduce the mechanism of Section \ref{Arnold},
we must find elementary transitions chains $a_0, \ldots, a_k$
in $\Rm^m$, with the additional requirement that $a_i\in A_{\mu}$.
It is necessary at this point to describe a bit more the set 
$A_{\mu}$.
Roughly, the KAM methods allow to prove the existence of the Torus $\mT_{\mu}(a)$
provided $a$ belongs to 
$$
A_{\mu}=\big\{a:\quad k\cdot (1,\omega_{\mu}(a))\geq \frac{\sqrt \mu}{\|k\|^{\tau}} \quad\forall k\in \Zm^{m+1}-\{0\}\big\}
$$
for some constant $\tau\geq m+1$.
This set $A_{\mu}$ is totally disconnected, hence it is not possible to
apply a continuity method like in   Proposition \ref{chain} in order to prove the
existence of a transition chain.
We must  be more quantitative, which is possible when $\mu$ is   so small that the Poincar\'e-Melnikov approximation is valid.
In that regime, we have $\bar \Delta_{a, \mu}\approx \mu M_a$, where $M_a$ has
 a non-degenerate critical point.
The conclusion of Lemma \ref{heteroclinics} can then be proved to hold 
 under the more explicit  condition that $\|a'-a\|\leq C\mu$.
In other words, the sequence $a_0, a_1, \ldots, a_k$ is an elementary transition chain
if $\|a_i-a_{i-1}\|\leq C\mu$.
However, the gaps in $A_{\mu}$ have a width of size $\sqrt{\mu}>C\mu$.
As a consequence, for small $\mu$, it seems impossible to build long transition chains,
and the method fails. 
This is the Large Gap Problem, 
see \cite{L}.
Even if there are classes of examples where the method can be applied because more tori
exist in some regions of phase space, see
\cite{BT,CG, BBB} for example, the generic case seems out of range.

\subsection{Normally hyperbolic invariant cylinder}\label{BAsec}
The Large Gap problem has now been solved, at least in the case where $m=1$,
see \cite{T, CY1, CY2, DLS, JAMS, Mo:02}.
We will not discuss and compare all these solutions here, but just
expose some general ideas which arise from them.

An important new point of view is to focus 
on the whole cylinder $\mC=\cup_a \mT (a)$ rather than on
each of the tori $\mT(a)$ individually.
This cylinder is Normally hyperbolic in the sense of \cite{HPS, F:71},
and thus it is preserved in the perturbed system. 
This new point of view is very natural, it appears in \cite{Mo:96,DLS:00}, 
and then in many other papers.
The deformed cylinder $\mC_{\mu}$ contains all the preserved tori
$\mT_{\mu}(a)$ obtained by KAM theory.
The restricted dynamics is described by an \textit{a priori}
stable system on $\Tm \times\Tm^m\times \Rm^m$.
If $m>1$, we are confronted to our lack of understanding of
the \textit{a priori} stable situation.
If $m=1$, however, the restricted system  is the suspension of an area preserving twist map, and we can exploit the good understanding of
these systems given by Birkhoff theory which has also
been interpreted (and extended) variationally in the works of Mather
 \cite{Ma:91, Ma:93}.
We consider this case ($m=1$) from now on.
The invariant $2$-tori which are graphs are of particular importance
(they correspond to rotational invariant circles of the time-one map).
To each of these invariant graphs, we can associate two real numbers,
the rotation number $\omega$ (defined from Poincar\'e theory 
of circle homeomorphisms), and the area $a$, which is the symplectic area
of the domain of the cylinder $\mC_{\mu}\cap\{t=0\}$ delimited by
the zero section and by the invariant graph under consideration.
If a given invariant graph $\mT$ of the restricted dynamics has irrational
rotation number (or is completely periodic), 
then there is no other invariant graph
with the same area $a$. We can take a 
two-covering and associate to this graph two functions $S_a^{\pm}$
by formula  similar to (\ref{Sa}).
 They generate the local  stable and unstable 
manifold of the Torus $\mT$, the correspond to the global
minima of the difference of the so-called barrier function  $S_a^+-S_a^-$.
Minimal homoclinics and isolated minimal homoclinics to $\mT$ can  be
defined as  in Section \ref{hd}. 
The existence of minimal homoclinics can be proved basically in the same way
as it was there.

\begin{defn}
An invariant graph is called a \textit{transition torus}
if it has irrational rotation number (or if it is foliated by
periodic orbits), and if all its minimal homoclinic orbits are isolated.
\end{defn}

Transition tori can be used to build transition chains
in the same way as partially hyperbolic  quasiperiodic tori with
transversal homoclinics.
Let $A\subset \Rm$ be the set of areas of transition  tori.
To each $a\in A$ is attached a unique transition torus $\mT_{\mu}(a)$
(note that this torus may be only Lipschitz, and is not necessarily
quasiperiodic).
If $A$ contains an interval $[a^-,a^+]$, then the existence
of a heteroclinic orbit between $\mT_{\mu}(a^-)$
and $\mT_{\mu}(a^+)$ can be proved by already exposed methods
(considering the way we have chosen our definitions, a variational method
should be used, but a parallel geometric theory could certainly be given).

In general, the set $A$ is totally disconnected, and transition chains
can't be obtained by a simple continuity method.
If we make the additional hypothesis that all invariant graphs of the 
restricted dynamics are transition tori, then the set $A$ is closed 
and a connected component $]a^-,a^+[$ of its complement corresponds
to a ``region of instability'' of the restricted system in the
 terminology of Birkhoff.
More precisely, the tori $\mT_{\mu}(a^-)$ and $\mT_{\mu}(a^+)$
enclose  a  cylinder which does not contain any invariant graph.
The theory of Birkhoff then implies that there exist orbits of the 
restricted dynamics connecting an arbitrarily small neighborhood
of $\mT_{\mu}(a^-)$ to an arbitrarily small neighborhood of $\mT_{\mu}(a^+)$.
This gives an indication about how to solve the large gap problem:
use the Birkhoff orbits to cross regions of instability, and the 
Arnold homoclinic mechanism to cross transition circles.
It is by no means obvious to prove the existence of actual orbits
 shadowing that kind of structure.
In order to do so, one should first put these mechanisms into
a common framework. The variational framework seems appropriate,
although a geometric approach is also possible.
The Birkhoff theory was described and extended using variational 
methods by Mather in \cite{Ma:91}, and he proposed a new variational
formalism adapted to higher dimensional situations in \cite{Ma:93}.
On the other hand, Bessi's method indicates how to put 
Arnold's mechanism into a variational framework. 
These heuristics lead to:

\begin{thm}\label{BA}
Let $[a^-,a^+]$ be a given interval.
If all the invariant graphs of area $a\in [a^-,a^+]$
of the restricted dynamics are transition tori,
then there exists an orbit
$(\theta(t), q(t), \dot \theta(t), \dot q(t))$ and a time $T>0$ such that 
$\dot \theta(0)\leq a^-$ and $\dot \theta (T)\geq a^+$.
\end{thm}

This theorem is proved using variational methods and weak KAM theory in 
\cite{JAMS}, Section 11, where it  is deduced from more general abstract results. It also almost follows from \cite{CY2}, Theorem 5.1,
which is another general abstract result proved by elaborations on
 Mather's variational methods \cite{Ma:93}, see also \cite{fourier}.
Applying that  result of Cheng and Yan, however, would
require  a minor  additional generic hypothesis
on the restricted dynamics. In the case where $r=1$, a slightly weaker
version of Theorem \ref{BA} could also be deduced from the  earlier paper of Chen and Yan \cite{CY1}. 
Under different sets of hypotheses, results in the same spirit have been 
obtained by geometric methods in \cite{GR:07,GR:09}.
At the moment, these methods do not reach statements as general as 
Theorem \ref{BA}, but they apply in contexts where the variational methods
 can't be used.

The following variant of Theorem \ref{BA} may deserve attention
in connection to the Arnold Mechanism:
Assume that $a^-$ and $a^+$ belong to $A$, or in other words
that there exist transition tori $\mT_{\mu}(a^{\pm})$.  
These tori enclose a compact invariant  piece $\mC_{\mu}[a^-,a^+] $ of the invariant
cylinder.
If all the invariant graphs contained in $\mC_{\mu}[a^-,a^+] $
are transition tori, then we say that $\mC_{\mu}[a^-,a^+] $ is
a transition channel.
The proof of Theorem \ref{BA} also implies that, if 
$\mC_{\mu}[a^-,a^+] $ is a transition channel, then there exists a 
heteroclinic orbit connecting $\mT_{\mu}(a^{-})$
to $\mT_{\mu}(a^{+})$.

Theorem \ref{BA} proves the existence of diffusion under ``explicit''
conditions. These conditions are hard to check on a given system,
but they seem to hold for typical systems.
It is much harder than one may expect  to prove a precise 
statement in that direction, but it was achieved by Cheng and Yan 
in \cite{CY1, CY2}.
The main difficulty comes from the condition on the isolated 
minimal homoclinics. Actually, it is not hard to prove that the homoclinics
to a given torus are isolated for a typical perturbation, but we need
the condition to hold for all the tori 
simultaneously.
Since there are uncountably many tori, it is necessary to understand the
regularity of the map  $a\lmto S^{\pm}_a$.
Recall that the functions $S^{\pm}_a$ are well-defined provided there
exists an invariant graph of area $a$ which has irrational rotation number 
or is foliated by periodic orbits. We call $\tilde A$ this set of areas,
it contains $A$.
Cheng and Yan prove that the map $a\lmto S_a^{\pm}$ is H\"older continuous on $\tilde A$,
and deduce the genericity result using an unpublished idea of John Mather.

\section{Back to the \textit{a priori} stable case}\label{hyper}

The main objects in Arnold's mechanism are partially hyperbolic 
tori, that he called whiskered tori.
It was proved by Treshchev \cite{T},
that whiskered tori exist in the \textit{a priori} stable 
situation, see also \cite{E:94,N}.
 However, because of the Large Gap Problem,
it seems difficult to prove directly the existence of transition chains
made of whiskered tori.
Actually, small transition chains do exist, because the density
of KAM tori increases near a given one, but the length of these chains
gets small when $\e$ gets small, hence these chains do not produce instability
of the action variables in general.

The modern paradigm on the \textit{a priori} unstable case that we exposed 
in Section \ref{BAsec}
 elects 3-dimensional normally hyperbolic invariant cylinders as
the important structure.
It is well-known that normally hyperbolic invariant cylinders exist
in the \textit{a priori} stable case.
For example, each $2$-dimensional whiskered torus has a center manifold,
which is a $3$-dimensional normally hyperbolic invariant cylinder, see 
\textit{e. g.}\cite{BT:00}.
Actually, it is simpler to prove directly the existence of normally
hyperbolic invariant cylinders, this involves no small
divisors. However, the most direct proofs seem to produce ``small''
normally hyperbolic cylinders, which means that their size is getting
small with $\e$, so that we  face the same problem 
as above when we had small transition chains.
The main statement of \cite{hyperbolic} is that ``large''
normally hyperbolic cylinders exist, meaning that their
size is bounded from below independently of $\e$.

In order to be more specific, let us select a resonant
 momentum of the form $y_0=(I_0,0)\in \Rm\times \Rm^{n-1}$,
with $I_0$ Diophantine. Assuming that the corresponding averaged
potential $V$ has a unique minimum and that this minimum is   non-degenerate, we have:

\begin{thm}[\cite{hyperbolic}]
There exists two intervals $[a^-,a^+]\subset J$, $J $ open, both 
\emph{independent from $\e$}, and $\epsilon_0>0$ such that, 
for $\e\in ]0,\e_0[$ the following holds:

There exists a $C^1$ map
$$(Q,P):\Tm\times \Tm \times J\ni (t,\theta, I)\lmto 
\big (Q(t,\theta,I),P(t,\theta, I)\big)\in \Tm^{n-1}\times \Rm^{n-1}
$$
such that the flow is tangent to the graph $\Gamma$ of $(Q,P)$.
Moreover, there exist two real numbers $a_0<a^-$  and $a_1>a^+$ in $J$
(which depend on $\e$)
such that the Treshchev tori $\mT(a_0)$ and $\mT(a_1)$ exist and are
contained in $\Gamma$.
The part $\Gamma_0^1$ of $\Gamma$ delimited by these two tori is then 
a compact invariant manifold with boundary of the flow, it is normally
hyperbolic. It is equivalent to say that it is partially hyperbolic
with a central distribution equal to the tangent space of $\Gamma$.
The inner dynamics is the suspension of an area-preserving twist map
(where the area is the one induced from the ambient symplectic form).
\end{thm}

It is then reasonable to expect that, under generic additional hypotheses,
$\Gamma_0^1$ is a transition channel as defined in Section \ref{BAsec},
and thus that $\mT(a_0)$ and $\mT(a_1)$ are connected by a heteroclinic 
orbit.
We are currently exploring that program in collaboration
with Vadim Kaloshin and Ke Zhang.
 It is important to observe that the map $(Q,P)$ is not $C^1$-close
to $(0,0)$, and that the inner dynamics is not close to integrable.
Fortunately,
Theorem \ref{BA} allows such a generality.

% Now the text of your article starts

% This is for the references %
% The reference below to Ahlfor's Complex Analysis book %
% is a sample that you should remove and put your own   %
% references %
% You can cite it with \cite{ahlfors} %

% Done :-) %
\end{document}